\documentclass[12pt,a4paper,reqno]{amsart}
\usepackage{amsmath, amsthm, amsfonts, amssymb, txfonts, graphicx, bm,url, verbatim,siunitx}

\theoremstyle{plain}
\newtheorem{thrm}{Theorem}[section]
\newtheorem{lmm}[thrm]{Lemma}
\newtheorem{prpstn}[thrm]{Proposition}
\newtheorem{crllry}[thrm]{Corollary}

\newtheorem*{rmk}{Remark}

\numberwithin{sblmm}{thrm} 
\numberwithin{equation}{section}

\usepackage[colorlinks=true, linkcolor=black, citecolor=black, urlcolor=black]{hyperref} 

\renewcommand{\phi}{\varphi}

\begin{document}
\title{On the Hardy--Ramanujan Theorem}
\author[B. Durkan]{Benjamin Durkan}
\address{Department of Mathematics, University of Manchester, Oxford Road, Manchester, M13 9PL}
\email{benjamin.durkan@manchester.ac.uk}

\begin{abstract}
    In this note we prove an effective version of the Hardy--Ramanujan Theorem. For every $x\ge 2$ and every non-negative function $F$ on the non-negative integers, we show 
    $$\frac{1}{x}\sum_{2\le n\le x}F(\omega(n)-1)\le 118\,\mathbb{E}F(Z_{\log\log x+4.096}),$$
    where $Z_{\lambda}$ is Poisson with parameter $\lambda$. Thus the shifted empirical distribution of $\omega(n)$ is pointwise dominated by a fixed multiple of a Poisson law. We also obtain the sharper squarefree analogue, derive explicit Chernoff and Gaussian-window estimates, obtain moderate-deviation upper bounds and uniform moment estimates, and transfer these consequences to $\Omega(n)$ and to the number of prime divisors occurring exactly once.
\end{abstract}

\maketitle
\section{Introduction}
For a positive integer $n$, let $$\omega(n)=\#\{p:p|n\}$$ denote the number of distinct prime factors of $n$, and let $$\Omega(n)=\sum_{p^a\lVert n}a$$ denote the total number of prime factors of $n$, counted with multiplicity. We also write
$$\omega_1(n)=\#\{p:p\lVert n\}.$$
Thus, if $n=\prod_{i=1}^rp_i^{a_i}$ with the $p_i$ distinct, then $\omega(n)=r$ and $\Omega(n)=a_1+\cdots+a_r$. It follows that $\omega(n)=\Omega(n)$ precisely when $n$ is squarefree. We use the convention $\omega(1)=\Omega(1)=\omega_1(1)=0$.

Hardy and Ramanujan \cite{HardyRamanujan1917} proved that both $\omega(n)$ and $\Omega(n)$ have normal order $\log\log n$. That is, given $\varepsilon>0$,
$$|\omega(n)-\log\log n|<\varepsilon\log\log n$$ for all but $o(x)$ integers $n\le x$.
Their proof in fact contains a sharper quantitative mechanism. It rests on a uniform counting inequality for the level sets of $\omega$. If
$$A_k(x):=\#\{n\le x:\omega(n)=k\},$$
then the classical Hardy--Ramanujan estimate has the form
\begin{equation*}
    A_k(x)\ll\frac{x}{\log x}\frac{(\log\log x+C_1)^{k-1}}{(k-1)!},
\end{equation*}
uniformly in $k$. This may be compared with Landau's fixed-$k$ theorem \cite{Landau1909},
$$A_k(x)\sim \frac{x}{\log x}\frac{(\log\log x)^{k-1}}{(k-1)!},$$
and with later refinements, notably the Sathe--Selberg theorem \cite{Sathe1953,Selberg1954}, which give asymptotics when $k$ is allowed to grow with $\log\log x$. See also Hildebrand and Tenenbaum \cite{HildebrandTenenbaum1988} for wider ranges and a detailed account of the local problem. This note has a different aim: it keeps the elementary Hardy--Ramanujan deletion argument, makes its constants explicit, and packages the conclusion as a domination inequality against a Poisson distribution.

This perspective places the original proof between two familiar successors. Tur\'{a}n's second-moment proof \cite{Turan1934} is shorter and already proves the normal order, but it gives only polynomial exceptional-set bounds by Chebyshev's inequality. The Erd\H{o}s--Kac theorem \cite{ErdosKac1940}, and later probabilistic treatments such as those of Kubilius \cite{Kubilius1964}, identify the limiting Gaussian law; see also Granville and Soundararajan \cite{GranvilleSoundararajan2007} for a modern sieve-oriented account. The Hardy--Ramanujan inequality does not prove a limiting distribution, but it gives exponential upper tails essentially for free once one recognizes the exponential-series structure.

Our main result is the following. Let $Z_{\lambda}$ denote a Poisson random variable with parameter $\lambda$.

\begin{thrm}\label{thm:effective_HR_ineq}
    Let $x\ge 2$, $\xi=\log\log x+4.096$ and $Q=1.95e^{4.096}<117.20$. Then for every non-negative function $F:\mathbb{Z}_{\ge 0}\to [0,\infty]$,
    \begin{equation}\label{eqn:thm_ineq}
        \frac{1}{x}\sum_{2\le n\le x}F(\omega(n)-1)\le Q\mathbb{E}F(Z_{\xi}),
    \end{equation}
    with the right-hand side interpreted in the extended sense.
\end{thrm}

Theorem \ref{thm:effective_HR_ineq} follows from the explicit counting estimate below.

\begin{thrm}\label{thm:HR_original_effective}
    For every real $x\ge 2$ and every integer $k\ge 1$,
    \begin{equation}\label{eqn:HR_original_effective}
        A_k(x)\le 1.95\frac{x}{\log x}\frac{(\log\log x+4.096)^{k-1}}{(k-1)!}.
    \end{equation}
\end{thrm}
The proof of Theorem \ref{thm:HR_original_effective} is elementary. We count the same integer several times by deleting exact prime-power divisors. This gives a recurrence for $A_k(x)$, and explicit estimates for prime sums supply the constants appearing in the theorem. We formulate this inductive approach as a general deletion principle. The shift $4.096$ is not meant to be best possible; its role is to keep the finite-$x$ argument transparent. What is more important is the functional form of Theorem \ref{thm:effective_HR_ineq}: since $F$ is arbitrary and non-negative, all pointwise, tail, generating-function and moment estimates follow by choosing $F$.

Several consequences are immediate. First, Theorem \ref{thm:effective_HR_ineq} gives large-deviation estimates, at fixed multiples of $\log\log x$, with the usual Poisson rate function.

\begin{crllry}\label{cor:fixed_mult_large_deviations}
    Let $I(\alpha)=\alpha\log\alpha-\alpha+1$ for $\alpha>0$. Then, for each fixed $\alpha>1$,
    \begin{equation}\label{eqn:cor_fixed_mult_1}
        \#\{2\le n\le x:\omega(n)\ge\alpha\log\log x\}\ll_{\alpha}x(\log x)^{-I(\alpha)}.
    \end{equation}
    For each fixed $0<\alpha<1$,
    \begin{equation}\label{eqn:cor_fixed_mult_2}
        \#\{2\le n\le x:\omega(n)\le\alpha\log\log x\}\ll_{\alpha}x(\log x)^{-I(\alpha)}.
    \end{equation}
\end{crllry}

Second, it gives Gaussian-scale exceptional-set bounds with explicit constants. These have the same shape as the Gaussian tail in the Erd\H{o}s--Kac theorem, although they are only upper bounds and do not imply convergence to the normal law.

\begin{thrm}\label{thm:effective_Gaussian}
    Suppose $x\ge e^{e^{200}}$. We have the following.
    \begin{enumerate}
        \item If $1\le A\le\sqrt{\log\log x}/4$, then
        \begin{align}
            \frac{1}{x}\#\{2\le n\le x:|&\omega(n)-\log\log x|>A\sqrt{\log\log x}\}\\&\le\frac{305}{A}\exp\left(-\frac{A^2}{2}+\frac{5.096A}{\sqrt{\log\log x}}+\frac{2.048A^2}{\log\log x}+\frac{A^3}{6\sqrt{\log\log x}}\right)\label{eqn:effec_gaussian_1}\\
            &\le\frac{305}{A}\exp\left(-\frac{A^2}{2}+\frac{7.31A^3}{\sqrt{\log\log x}}\right).\notag
        \end{align}
        \item If $1\le A\le\sqrt{\log\log x}/12$, then
        \begin{equation}\label{eqn:effec_gaussian_2}
            \frac{1}{x}\#\{2\le n\le x:|\Omega(n)-\log\log x|>A\sqrt{\log\log x}\}\le\frac{332}{A}\exp\left(-\frac{A^2}{2}+\frac{2.29A^3}{\sqrt{\log\log x}}\right).
        \end{equation}
        The same estimate holds with $\Omega$ replaced by $\omega_1$.
    \end{enumerate}
\end{thrm}
The preceding explicit estimates imply the following more standard moderate-deviation formulation.
\begin{thrm}\label{thm:mod_dev_upper_bound}
    Suppose $A=A(x)\to\infty$ satisfies $A=o(\sqrt{\log\log x})$. Then for $f=\omega$, $f=\Omega$ or $f=\omega_1$,
    \begin{equation}\label{eqn:mod_dev_upper_bound}
        \frac{1}{x}\#\{2\le n\le x:|f(n)-\log\log x|>A\sqrt{\log\log x}\}\le\exp\left(-\left(\frac{1}{2}+o(1)\right)A^2\right).
    \end{equation}
\end{thrm}
The domination theorem also gives useful uniform moment estimates.

\begin{thrm}\label{thm:uniform_moment_bounds}
    There is an absolute constant $C>0$ such that, for every $x\ge 3$, every real $r\ge 1$ and $f=\omega$, $f=\Omega$ or $f=\omega_1$ we have
    \begin{equation}\label{eqn:uniform_moments}
        \left(\frac{1}{x}\sum_{2\le n\le x}|f(n)-\log\log x|^r\right)^{1/r}\le C(\sqrt{r(\log\log x+1)}+r).
    \end{equation}
\end{thrm}
For fixed $r$, Theorem \ref{thm:uniform_moment_bounds} recovers the Gaussian-order central moment estimate
\begin{equation}\label{eqn:gaussian_order_central_moment}
    \frac{1}{x}\sum_{2\le n\le x}|f(n)-\log\log x|^r\ll_r(\log\log x+1)^{r/2}.
\end{equation}
The normal-order theorem follows from Theorem \ref{thm:effective_Gaussian}, or more softly from \eqref{eqn:gaussian_order_central_moment} together with Markov's inequality.

\begin{thrm}\label{thm:original_HR_theorem}
    Let $\psi(n)$ be any function tending to infinity. Then, for all but $o(x)$ integers $n\le x$,
    $$|\omega(n)-\log\log n|\le\psi(n)\sqrt{\log\log n}.$$
    The same assertion holds with $\omega(n)$ replaced by $\Omega(n)$ or by $\omega_1(n)$.
\end{thrm}
Tur\'{a}n \cite{Turan1934} later gave a shorter proof of the normal order theorem by establishing the variance bound
$$\frac{1}{x}\sum_{n\le x}(\omega(n)-\log\log x)^2\ll \log\log x,$$
and then appealing to Chebyshev's inequality. Tur\'{a}n's proof was far shorter than Hardy and Ramanujan's, but has the drawback of giving only the polynomial tail 
$$\frac{1}{x}\#\{n\le x:|\omega(n)-\log\log x|>A\sqrt{\log\log x}\}\ll\frac{1}{A^2}.$$
The original Hardy--Ramanujan counting inequality contains sharper information.

\section{Preliminary effective estimates}
Throughout the paper, we shall use $p$ to denote a prime number and $\log$ to denote the natural logarithm. We shall use the following explicit estimates from prime number theory.

\begin{lmm}\label{lem:effective_prime_estimates}
    For all $x\ge 2$ we have the following:
    \begin{align}
        &\sum_{p\le x}\frac{1}{p}\le\log\log x+0.867,\label{eqn:mertens_explicit}\\
        &\sum_{p\le x}\frac{\log p}{p}<\log x,\label{eqn:chebyshev_explicit}\\
       &\pi(x)<1.25506\frac{x}{\log x}.\label{eqn:pi_explicit}
    \end{align}
    Moreover
    $$J:=\sum_{a=2}^{\infty}\sum_p\frac{a+1}{p^a}<2.922.$$
\end{lmm}
\begin{proof}
    Rosser and Schoenfeld \cite{RosserSchoenfeld1962} prove the middle two statements, and also
    $$\sum_{p\le x}\frac{1}{p}\le\log\log x+B_1+\frac{1}{2\log^2x}\qquad (x\ge 286),$$
    where $B_1=0.26149\cdots$ is the Meissel--Mertens constant. This is stronger than the first estimate for $x\ge 286$. For $2\le x<286$, the function $\sum_{p\le x}\frac{1}{p}-\log\log x$ is decreasing between consecutive primes and therefore needs only be checked at primes below $286$. The maximum occurs at $x=2$, where it equals $1/2-\log\log 2=0.8665\cdots$, establishing the first bound. For the final bound, we note that
    $$\sum_{a=2}^{\infty}\frac{a+1}{p^a}=\frac{3p-2}{p(p-1)^2}.$$
    We compute
    $$\sum_{p\le 5000}\frac{3p-2}{p(p-1)^2}<2.921316,$$
    and for $p>5000$ we have
    $$\frac{3p-2}{p(p-1)^2}<\frac{3}{(p-1)^2},$$
    and so the remaining tail is at most
    $$\sum_{n\ge 5001}\frac{3}{(n-1)^2}<\frac{3}{4999}<0.000601,$$ so that $J<2.922$.
\end{proof}
\section{The Hardy--Ramanujan counting inequality}
For $k\ge 1$, write
$$A_k(x):=\#\{n\le x:\omega(n)=k\}.$$
The proof of the counting inequality is based on the following deletion recurrence.

\begin{lmm}\label{lem:deletion_recurrence}
    For every integer $\nu\ge 1$ and every real $x\ge 2$,
    \begin{equation}\label{eqn:deletion_recurrence}
    \nu A_{\nu+1}(x)\le\sum_{\substack{a\ge 1\\ p^{a+1}\le x}}A_{\nu}\left(\frac{x}{p^a}\right).
    \end{equation}
\end{lmm}
\begin{proof}
    Let $n=q_1^{b_1}q_2^{b_2}\cdots q_{\nu+1}^{b_{\nu+1}}\le x$, where $q_1<\cdots<q_{\nu+1}$ are primes. Then $\omega(n)=\nu+1$. For each $1\le i\le\nu$, delete the exact power $q_i^{b_i}$. The remaining integer has exactly $\nu$ distinct prime factors. Moreover, $q_i^{b_i+1}\le q_i^{b_i}q_{\nu+1}\le n\le x$, so the pair $(q_i,b_i)$ is admissible in the sum on the right-hand side of \eqref{eqn:deletion_recurrence}. Thus each $n$ counted by $A_{\nu+1}(x)$ contributes at least $\nu$ admissible deletions.
\end{proof}

The recurrence may be recast as a general induction principle.

\begin{prpstn}\label{prop:deletion_machine}
    Let $C>0$ and $B\in\mathbb{R}$. Suppose that $\log\log y+B>0$ for every $y\ge 2$, and that the following two estimates hold for every $y\ge 2$:
    \begin{equation}\label{eqn:del_rec_1}
        A_1(y)\le C\frac{y}{\log y},
    \end{equation}
    \begin{equation}\label{eqn:del_rec_2}
        \sum_{\substack{a\ge 1\\ p^{a+1}\le y}}\frac{1}{p^a\log(y/p^a)}\le\frac{\log\log y+B}{\log y}.
    \end{equation}
    Then, for every $y\ge 2$ and every integer $k\ge 1$,
    $$A_k(y)\le\frac{Cy}{\log y}\frac{(\log\log y+B)^{k-1}}{(k-1)!}.$$
\end{prpstn}
\begin{proof}
    The case $k=1$ is \eqref{eqn:del_rec_1}. Suppose that the proposition has been proved for $k=\nu$, uniformly for all $y\ge 2$. By Lemma \ref{lem:deletion_recurrence} we have
    $$\nu A_{\nu+1}(x)\le\sum_{\substack{a\ge 1\\ p^{a+1}\le x}}A_{\nu}\left(\frac{x}{p^a}\right).$$
    Since $p^{a+1}\le x$ implies $x/p^a\ge p\ge 2$, the induction hypothesis gives
    $$A_{\nu}\left(\frac{x}{p^a}\right)\le C\frac{x/p^a}{\log(x/p^a)}\frac{(\log\log(x/p^a)+B)^{\nu-1}}{(\nu-1)!}.$$
    Moreover, since $2\le x/p^a\le x$, the factor $\log\log(x/p^a)+B$ is at most $\log\log x+B$. It follows that
    $$A_{\nu+1}(x)\le Cx\frac{(\log\log x+B)^{\nu-1}}{\nu!}\sum_{\substack{a\ge 1\\ p^{a+1}\le x}}\frac{1}{p^a\log(x/p^a)}.$$
    Substituting in \eqref{eqn:del_rec_2} completes the proof.
\end{proof}
It remains to verify the two hypotheses of Proposition \ref{prop:deletion_machine} with explicit constants.

\begin{lmm}\label{lem:effective_1}
  For every $x\ge 2$,
  $$\sum_{\substack{a\ge 1\\ p^{a+1}\le x}}\frac{1}{p^a\log(x/p^a)}\le\frac{\log\log x+4.096}{\log x}.$$
\end{lmm}
\begin{proof}
    Since the sum is empty for $x<4$ we may assume that $x\ge 4$. For $a=1$ the condition $p^2\le x$ gives $\log p\le\frac{1}{2}\log x$. Hence
    $$\frac{1}{\log(x/p)}=\frac{1}{\log x-\log p}\le\frac{1}{\log x}+\frac{2\log p}{(\log x)^2}.$$
    Applying Lemma \ref{lem:effective_prime_estimates} gives
    \begin{align*}
        \sum_{p^2\le x}\frac{1}{p\log(x/p)}&\le\frac{\log\log\sqrt{x}+0.867}{\log x}+\frac{2}{(\log x)^2}\sum_{p\le\sqrt{x}}\frac{\log p}{p}\\
        &<\frac{\log\log x-\log 2+0.867}{\log x}+\frac{1}{\log x}\\
        &<\frac{\log\log x+1.174}{\log x}.
    \end{align*}
    For $a\ge 2$, the condition $p^{a+1}\le x$ implies $\log(x/p^a)\ge\frac{\log x}{a+1}$, so by Lemma \ref{lem:effective_prime_estimates} we have
    $$\sum_{\substack{a\ge 2\\ p^{a+1}\le x}}\frac{1}{p^a\log(x/p^a)}\le\frac{1}{\log x}\sum_{a=2}^{\infty}\sum_p\frac{a+1}{p^a}<\frac{2.922}{\log x}.$$
    Since $1.174+2.922<4.096$ the result follows.
\end{proof}

\begin{lmm}\label{lem:prime_powers}
   For every real $x\ge 2$,
   $$A_1(x)\le 1.95\frac{x}{\log x}.$$
\end{lmm}
\begin{proof}
    The integers counted by $A_1(x)$ are precisely the prime powers $p^a\le x$, each with a unique representation. Thus $A_1(x)=\sum_{a=1}^{\infty}\pi(x^{1/a})$, where only terms with $a\le\log x/\log 2$ contribute. For $x\ge e^{10}$, Lemma \ref{lem:effective_prime_estimates} yields
    $$\pi(x)+\pi(\sqrt{x})<1.25506\frac{x}{\log x}+2.51012\frac{\sqrt{x}}{\log x}.$$
    Also
    $$\sum_{3\le a\le \log x/\log 2}x^{1/a}\le x^{1/3}+\frac{\log x}{\log 2}x^{1/4}.$$
    After division by $x/\log x$, the last three contributions are bounded by
    $$2.51012x^{-1/2}+(\log x)x^{-2/3}+\frac{(\log x)^2}{\log 2}x^{-3/4},$$
    which is decreasing for $x\ge e^{10}$ and is $<0.11$ at $x=e^{10}$. Thus $A_1(x)<1.37\frac{x}{\log x}$ for $x\ge e^{10}$. It remains to check $2\le x<e^{10}$. Between consecutive prime powers, $A_1(x)$ is constant, while $\log x/x$ is increasing on $(1,e)$ and decreasing on $(e,\infty)$. A finite check at the relevant endpoints shows that the maximum of $A_1(x)\log x/x$ occurs at $x=32$, where
    $$A_1(32)\frac{\log 32}{32}=18\frac{\log 32}{32}=1.949476\cdots<1.95.$$
    This proves Lemma \ref{lem:prime_powers}.
\end{proof}

\begin{proof}[Proof of Theorem \ref{thm:HR_original_effective}]
    Apply Proposition \ref{prop:deletion_machine} with $C=1.95$ and $B=4.096$, using Lemmas \ref{lem:effective_1} and \ref{lem:prime_powers}.
\end{proof}

\begin{rmk}
    Hardy and Ramanujan first treated squarefree integers and then passed to general integers. The proof above goes directly to unrestricted integers by counting exact prime-power divisors. This keeps the original deletion idea while avoiding a separate squarefree argument.
\end{rmk}

The same argument gives a sharper estimate on the squarefree integers. Write
$$S_k(x):=\#\{n\le x:\mu^2(n)=1,\ \omega(n)=k\}.$$

\begin{prpstn}\label{prop:squarefree_counting}
    For every real $x\ge 2$ and every integer $k\ge 1$,
    \begin{equation}\label{eqn:squarefree_counting}
        S_k(x)\le 1.25506\frac{x}{\log x}\frac{(\log\log x+1.174)^{k-1}}{(k-1)!}.
    \end{equation}
\end{prpstn}
\begin{proof}
    If $n=q_1\cdots q_{\nu+1}\le x$ is squarefree with $q_1<\cdots<q_{\nu+1}$, then deleting any of $q_1,\ldots,q_{\nu}$ leaves a squarefree integer with $\nu$ prime factors, and $q_i^2\le q_iq_{\nu+1}\le n\le x$. Thus
    $$\nu S_{\nu+1}(x)\le\sum_{p^2\le x}S_{\nu}\left(\frac{x}{p}\right).$$
    The same induction as in Proposition \ref{prop:deletion_machine}, with $C=1.25506$ and with only the $a=1$ part of Lemma \ref{lem:effective_1}, applies. The base case is $S_1(x)=\pi(x)<1.25506x/\log x$, and the needed sum is
    $$\sum_{p^2\le x}\frac{1}{p\log(x/p)}\le\frac{\log\log x+1.174}{\log x},$$
    which is the $a=1$ estimate in the proof of Lemma \ref{lem:effective_1}.
\end{proof}

\begin{crllry}\label{cor:squarefree_domination}
    Let $x\ge 2$, $\xi_0=\log\log x+1.174$ and $Q_0=1.25506e^{1.174}<4.07$. Then, for every non-negative function $F:\mathbb{Z}_{\ge 0}\to[0,\infty]$,
    \begin{equation}\label{eqn:squarefree_domination}
        \frac{1}{x}\sum_{\substack{2\le n\le x\\ \mu^2(n)=1}}F(\omega(n)-1)\le Q_0\mathbb{E}F(Z_{\xi_0}).
    \end{equation}
\end{crllry}
\begin{proof}
    This is obtained from Proposition \ref{prop:squarefree_counting} exactly as Theorem \ref{thm:effective_HR_ineq} is obtained from Theorem \ref{thm:HR_original_effective}.
\end{proof}

\section{Poisson domination and exponential tails}
Let $Z_{\lambda}$ denote a Poisson random variable of mean $\lambda$, so that 
$$\mathbb{P}(Z_{\lambda}=m)=e^{-\lambda}\frac{\lambda^m}{m!} \qquad (m=0,1,2,\cdots).$$

Throughout this section we set $L=\log\log x$, $\xi=L+4.096$ and $Q=1.95e^{4.096}<117.20$.
\begin{proof}[Proof of Theorem \ref{thm:effective_HR_ineq}]
    For $m\ge 0$, Theorem \ref{thm:HR_original_effective} gives
    $$A_{m+1}(x)\le 1.95\frac{x}{\log x}\frac{\xi^m}{m!}.$$

    Since $e^{\xi}=e^{4.096}\log x$, we may rewrite this as
    $$\frac{A_{m+1}(x)}{x}\le 1.95e^{4.096}e^{-\xi}\frac{\xi^m}{m!}=Q\mathbb{P}(Z_{\xi}=m).$$
    Therefore, for every non-negative $F:\mathbb{Z}_{\ge 0}\to [0,\infty]$, we have
    \begin{align*}
        \frac{1}{x}\sum_{2\le n\le x}F(\omega(n)-1)&=\frac{1}{x}\sum_{m=0}^{\infty}F(m)A_{m+1}(x)\\
        &\le Q\sum_{m=0}^{\infty}F(m)e^{-\xi}\frac{\xi^m}{m!}\\
        &=Q\mathbb{E}F(Z_{\xi}),
    \end{align*}
    as claimed.
\end{proof}

Taking $F=\mathbf{1}_S$ gives the setwise form
\begin{equation}\label{eqn:setwise}
 \#\{2\le n\le x:\omega(n)-1\in S\}\le Qx\mathbb{P}(Z_{\xi}\in S)
\end{equation}
for every $S\subseteq\mathbb{Z}_{\ge 0}$. Taking $F(m)=e^{\theta m}$ gives the effective moment-generating-function bound
\begin{equation}\label{eqn:mgf}
    \frac{1}{x}\sum_{2\le n\le x}e^{\theta(\omega(n)-1)}\le Q\exp\left(\xi(e^{\theta}-1)\right) \qquad (\theta\in\mathbb{R}).
\end{equation}
We now define two functions which we shall use repeatedly. Define, for $u\ge 0$,
$$h_1(u)=(1+u)\log(1+u)-u$$
and, for $0\le u<1$,
$$h_2(u)=(1-u)\log(1-u)+u.$$

\begin{prpstn}\label{prop:chernoff}
    For every $t\ge 0$,
    $$\frac{1}{x}\#\{2\le n\le x:\omega(n)-1-\xi\ge t\}\le Q\exp\left(-\xi h_1\left(\frac{t}{\xi}\right)\right).$$
    For every $0\le t<\xi$,
    $$\frac{1}{x}\#\{2\le n\le x:\omega(n)-1-\xi\le -t\}\le Q\exp\left(-\xi h_2\left(\frac{t}{\xi}\right)\right).$$
\end{prpstn}
\begin{proof}
    For the upper tail, Markov's inequality and \eqref{eqn:mgf} give, for $\theta>0$,
    $$\frac{1}{x}\#\{2\le n\le x:\omega(n)-1\ge\xi+t\}\le Q\exp\left(\xi(e^{\theta}-1)-\theta(\xi+t)\right).$$
    Optimising with $e^{\theta}=1+t/\xi$ gives the first bound. The lower tail is identical, using $e^{-\theta(\omega(n)-1)}$ and optimising with $e^{-\theta}=1-t/\xi$.
\end{proof}
The Chernoff bounds may also be written in the form of Bernstein inequalities.

\begin{crllry}\label{cor:bernstein}
    For every $t\ge 0$,
    $$\frac{1}{x}\#\{2\le n\le x:\omega(n)-1-\xi\ge t\}\le Q\exp\left(-\frac{t^2}{2(\xi+t/3)}\right)$$
    and for every $0\le t<\xi$,
    $$\frac{1}{x}\#\{2\le n\le x:\omega(n)-1-\xi\le -t\}\le Q\exp\left(-\frac{t^2}{2\xi}\right).$$
\end{crllry}
\begin{proof}
    Use the inequalities $h_1(u)\ge\frac{u^2}{2(1+u/3)}$ for $u\ge 0$ and $h_2(u)\ge\frac{u^2}{2}$ for $0\le u<1$, in Proposition \ref{prop:chernoff}.
\end{proof}

\begin{proof}[Proof of Corollary \ref{cor:fixed_mult_large_deviations}]
    Let $\alpha>1$ be fixed. For all sufficiently large $x$, the event $\omega(n)\ge\alpha\log\log x$ implies $\omega(n)-1-\xi\ge(\alpha-1)\log\log x-5.096$. Proposition \ref{prop:chernoff} gives an upper bound of the form
    $$\exp(-L I(\alpha)+O_{\alpha}(1)),$$
    because
    $$\xi h_1\left(\frac{(\alpha-1)L-5.096}{\xi}\right)=L(\alpha\log\alpha-\alpha+1)+O_{\alpha}(1).$$
    This proves \eqref{eqn:cor_fixed_mult_1}. The case $0<\alpha<1$ follows in the same way from Proposition \ref{prop:chernoff}, using
    $$\xi h_2\left(\frac{(1-\alpha)L+5.096}{\xi}\right)=L(\alpha\log\alpha-\alpha+1)+O_{\alpha}(1).$$
    The finitely many remaining values of $x$ are absorbed by the implicit constants.
\end{proof}

For comparison with Gaussian tails, we record the following estimate for a Poisson random variable.

\begin{lmm}\label{lem:poisson_mills}
    Let $Z_{\lambda}$ be Poisson with mean $\lambda$. If $\lambda\ge 16$ and $0<B\le\sqrt{\lambda}/2$, then
    \begin{equation}\label{eqn:poisson_upper_gaussian}
        \mathbb{P}(Z_{\lambda}-\lambda\ge B\sqrt{\lambda})\le\frac{1}{B}\exp\left(-\frac{B^2}{2}+\frac{B^3}{6\sqrt{\lambda}}\right),
    \end{equation}
    and
    \begin{equation}\label{eqn:poisson_lower_gaussian}
        \mathbb{P}(Z_{\lambda}-\lambda\le -B\sqrt{\lambda})\le\frac{1}{B}\exp\left(-\frac{B^2}{2}\right).
    \end{equation}
\end{lmm}
\begin{proof}
    Put $h(u)=(1+u)\log(1+u)-u$. For the upper tail, let
    $$m=\left\lceil \lambda+B\sqrt{\lambda}\right\rceil.$$
    For $k\ge m$,
    $$\frac{\mathbb{P}(Z_{\lambda}=k+1)}{\mathbb{P}(Z_{\lambda}=k)}=\frac{\lambda}{k+1}\le\frac{\lambda}{m+1}.$$
    Thus
    $$\mathbb{P}(Z_{\lambda}\ge m)\le\frac{\mathbb{P}(Z_{\lambda}=m)}{1-\lambda/(m+1)}.$$
    Since $B\le\sqrt{\lambda}/2$ and $\lambda\ge 16$, we have
    $$m+1\le\lambda+B\sqrt{\lambda}+2\le 2\lambda,\qquad m+1-\lambda\ge B\sqrt{\lambda},$$
    and hence
    \begin{equation}\label{eqn:geom_upper_gap}
        1-\frac{\lambda}{m+1}\ge\frac{B}{2\sqrt{\lambda}}.
    \end{equation}
    Stirling's lower bound $m!\ge\sqrt{2\pi m}(m/e)^m$, together with $m\ge\lambda$, gives
    $$\mathbb{P}(Z_{\lambda}=m)\le\frac{1}{\sqrt{2\pi m}}\exp\left(-\lambda h\left(\frac{m-\lambda}{\lambda}\right)\right).$$
    Combining this with \eqref{eqn:geom_upper_gap} yields
    $$\mathbb{P}(Z_{\lambda}\ge m)\le\frac{1}{B}\exp\left(-\lambda h\left(\frac{m-\lambda}{\lambda}\right)\right),$$
    because $2\sqrt{\lambda}/\sqrt{2\pi m}<1$. Since $h'(u)=\log(1+u)\ge 0$ and $h(u)\ge u^2/2-u^3/6$ for $u\ge 0$,
    $$h\left(\frac{m-\lambda}{\lambda}\right)\ge h\left(\frac{B}{\sqrt{\lambda}}\right)\ge\frac{B^2}{2\lambda}-\frac{B^3}{6\lambda^{3/2}},$$
    which proves \eqref{eqn:poisson_upper_gaussian}.

    For the lower tail, let
    $$m=\left\lfloor \lambda-B\sqrt{\lambda}\right\rfloor.$$
    For $k\le m$,
    $$\frac{\mathbb{P}(Z_{\lambda}=k-1)}{\mathbb{P}(Z_{\lambda}=k)}=\frac{k}{\lambda}\le\frac{m}{\lambda},$$
    and hence
    $$\mathbb{P}(Z_{\lambda}\le m)\le\frac{\mathbb{P}(Z_{\lambda}=m)}{1-m/\lambda}.$$
    Here $1-m/\lambda\ge B/\sqrt{\lambda}$. Also $m\ge\lambda/2-1\ge\lambda/4$, since $\lambda\ge 16$. Stirling's bound gives
    $$\mathbb{P}(Z_{\lambda}=m)\le\frac{1}{\sqrt{2\pi m}}\exp(-\lambda h_2(s)),\qquad s=\frac{\lambda-m}{\lambda}.$$
    Thus
    $$\mathbb{P}(Z_{\lambda}\le m)\le\frac{1}{B}\exp(-\lambda h_2(s)).$$
    The function $h_2$ is increasing on $[0,1)$, and $h_2(u)\ge u^2/2$ for $0\le u<1$. Since $s\ge B/\sqrt{\lambda}$, \eqref{eqn:poisson_lower_gaussian} follows.
\end{proof}

\begin{prpstn}\label{prop:gaussian_omega}
    Let $L=\log\log x$. If $L\ge 200$ and $1\le A\le\sqrt{L}/4$, then
    \begin{align}
        \frac{1}{x}\#\{2\le n\le x:\omega(n)-L>A\sqrt{L}\}&\le\frac{8Q}{5A}\exp\left(-\frac{A^2}{2}+\frac{5.096A}{\sqrt{L}}+\frac{2.048A^2}{L}+\frac{A^3}{6\sqrt{L}}\right),\label{eqn:omega_upper_window}\\
        \frac{1}{x}\#\{2\le n\le x:L-\omega(n)>A\sqrt{L}\}&\le\frac{Q}{A}e^{-A^2/2}.\label{eqn:omega_lower_window}
    \end{align}
    Consequently,
    \begin{multline}\label{eqn:omega_two_sided_window}
        \frac{1}{x}\#\{2\le n\le x:|\omega(n)-L|>A\sqrt{L}\}\\
        \le\frac{305}{A}\exp\left(-\frac{A^2}{2}+\frac{5.096A}{\sqrt{L}}+\frac{2.048A^2}{L}+\frac{A^3}{6\sqrt{L}}\right)\\
        \le\frac{305}{A}\exp\left(-\frac{A^2}{2}+\frac{7.31A^3}{\sqrt{L}}\right).
    \end{multline}
\end{prpstn}
\begin{proof}
    Recall that $\xi=L+4.096$, and put $d=5.096$. If $\omega(n)-L>A\sqrt{L}$, then
    $$(\omega(n)-1)-\xi>A\sqrt{L}-d.$$
    Set
    $$B_+=\frac{A\sqrt{L}-d}{\sqrt{\xi}}.$$
    The assumptions imply $B_+>0$, and in fact $B_+\ge 5A/8$. Indeed, after division by $A$, the left side of the desired inequality is minimized at $A=1$, and
    $$\sqrt{L}-\frac{5}{8}\sqrt{L+4.096}$$
    is increasing for $L\ge 200$ and is already $>5.096$ at $L=200$. Also
    $$B_+\le\frac{A\sqrt{L}}{\sqrt{\xi}}\le\frac{\sqrt{\xi}}{2}.$$
    The setwise domination \eqref{eqn:setwise} and Lemma \ref{lem:poisson_mills} give
    $$\frac{1}{x}\#\{2\le n\le x:\omega(n)-L>A\sqrt{L}\}\le\frac{8Q}{5A}\exp\left(-\frac{B_+^2}{2}+\frac{B_+^3}{6\sqrt{\xi}}\right).$$
    Now
    $$-\frac{B_+^2}{2}=-\frac{(A\sqrt{L}-d)^2}{2(L+4.096)}\le -\frac{A^2}{2}+\frac{dA}{\sqrt{L}}+\frac{2.048A^2}{L},$$
    and
    $$\frac{B_+^3}{6\sqrt{\xi}}\le\frac{A^3}{6\sqrt{L}}.$$
    This proves \eqref{eqn:omega_upper_window}.

    For the lower tail, if $L-\omega(n)>A\sqrt{L}$, then
    $$(\omega(n)-1)-\xi<-(A\sqrt{L}+d).$$
    Set
    $$B_- =\frac{A\sqrt{L}+d}{\sqrt{\xi}}.$$
    We have $B_-\ge A$, since
    $$A(\sqrt{\xi}-\sqrt{L})=\frac{4.096A}{\sqrt{\xi}+\sqrt{L}}\le\frac{4.096\sqrt{L}}{4(\sqrt{\xi}+\sqrt{L})}<0.52<d.$$
    Also $B_-\le\sqrt{\xi}/2$, since
    $$A\sqrt{L}+d\le\frac{L}{4}+5.096\le\frac{L+4.096}{2}$$
    for $L\ge 200$. Lemma \ref{lem:poisson_mills} gives \eqref{eqn:omega_lower_window}. Adding the two one-sided estimates gives the first inequality in \eqref{eqn:omega_two_sided_window}, since $(8/5+1)Q<305$. The simplified estimate follows from
    $$\frac{5.096A}{\sqrt{L}}+\frac{2.048A^2}{L}+\frac{A^3}{6\sqrt{L}}\le\frac{7.31A^3}{\sqrt{L}}$$
    for $A\ge 1$ and $L\ge 1$.
\end{proof}

\begin{proof}[Proof of Theorem \ref{thm:effective_Gaussian} for $\omega$]
    This is Proposition \ref{prop:gaussian_omega}.
\end{proof}

\section{The squarefull excess and the functions \texorpdfstring{$\Omega$ and $\omega_1$}{Omega and omega1}}
The passage from $\omega$ to $\Omega$ is controlled by
$$R(n)=\Omega(n)-\omega(n)=\sum_p(v_p(n)-1)_+.$$
We shall use an effective exponential tail for $R(n)$, rather than only its bounded mean.

\begin{prpstn}\label{prop:R_tail}
    For all $x\ge 1$ and all integers $T\ge 0$,
    \begin{equation}\label{eqn:R_tail}
        \frac{1}{x}\#\{n\le x:R(n)\ge T\}\le\frac{11}{2}2^{-T}.
    \end{equation}
\end{prpstn}
\begin{proof}
    Let $\mathcal{F}$ be the set of powerful numbers. If $R(n)\ge T$, then the powerful part
    $$d=\prod_{\substack{p^a\lVert n\\ a\ge 2}}p^a$$
    satisfies $d|n$, $d\in\mathcal{F}$ and $R(d)=R(n)\ge T$. Therefore
    $$\#\{n\le x:R(n)\ge T\}\le\sum_{\substack{d\in\mathcal{F}\\ R(d)\ge T}}\left\lfloor\frac{x}{d}\right\rfloor\le x\sum_{\substack{d\in\mathcal{F}\\ R(d)\ge T}}\frac{1}{d}.$$
    It remains to bound the reciprocal sum.

    Let $\mathcal{F}_{\rm odd}$ be the set of odd powerful numbers and put
    $$H=\sum_{e\in\mathcal{F}_{\rm odd}}\frac{2^{R(e)}}{e}.$$
    By multiplicativity,
    $$H=\prod_{p>2}\left(1+\sum_{a\ge 2}\frac{2^{a-1}}{p^a}\right)=\prod_{p>2}\left(1+\frac{2}{p(p-2)}\right).$$
    Multiplying the factors with $3\le p\le 97$, and bounding the rest by the corresponding integer product, gives
    $$H\le\prod_{3\le p\le 97}\left(1+\frac{2}{p(p-2)}\right)\exp\left(\sum_{m\ge 98}\frac{2}{m(m-2)}\right)<2.2.$$

    Fix $e\in\mathcal{F}_{\rm odd}$ and write $r=R(e)$. Every powerful number can be written uniquely as $2^ae$, with $e\in\mathcal{F}_{\rm odd}$ and either $a=0$ or $a\ge 2$. If $r\ge T$, the allowed powers of $2$ have total reciprocal weight
    $$1+\sum_{a\ge 2}2^{-a}=\frac{3}{2}.$$
    If $r<T$, then we need $a-1\ge T-r$, and the allowed powers of $2$ have reciprocal weight
    $$\sum_{a\ge T-r+1}2^{-a}=2^{-(T-r)}.$$
    Hence
    \begin{align*}
        \sum_{\substack{d\in\mathcal{F}\\ R(d)\ge T}}\frac{1}{d}&\le\frac{3}{2}\sum_{\substack{e\in\mathcal{F}_{\rm odd}\\ R(e)\ge T}}\frac{1}{e}+\sum_{\substack{e\in\mathcal{F}_{\rm odd}\\ R(e)<T}}\frac{2^{-(T-R(e))}}{e}\\
        &\le\frac{5}{2}H2^{-T}<\frac{11}{2}2^{-T}.
    \end{align*}
    This proves \eqref{eqn:R_tail}.
\end{proof}

\begin{prpstn}\label{prop:gaussian_Omega}
    If $L=\log\log x$, $L\ge 200$ and $1\le A\le\sqrt{L}/12$, then
    \begin{equation}\label{eqn:Omega_window}
        \frac{1}{x}\#\{2\le n\le x:|\Omega(n)-L|>A\sqrt{L}\}\le\frac{332}{A}\exp\left(-\frac{A^2}{2}+\frac{2.29A^3}{\sqrt{L}}\right).
    \end{equation}
\end{prpstn}
\begin{proof}
    For $1\le A\le 3$, the right side of \eqref{eqn:Omega_window} is $>1$, so the estimate is trivial. Assume $A>3$. Define
    $$T_A=\left\lceil\frac{A^2/2+\log(11eA/2)}{\log 2}\right\rceil.$$
    By Proposition \ref{prop:R_tail},
    \begin{equation}\label{eqn:R_tail_TA}
        \frac{1}{x}\#\{n\le x:R(n)\ge T_A\}\le\frac{1}{eA}e^{-A^2/2}.
    \end{equation}
    We shall use
    \begin{equation}\label{eqn:TA_bound}
        T_A\le\frac{3}{2}A^2\qquad (A\ge 3).
    \end{equation}
    Indeed, after accounting for the ceiling, \eqref{eqn:TA_bound} follows from
    $$\log(11eA/2)+\log 2\le\left(\frac{3}{2}\log 2-\frac{1}{2}\right)A^2,$$
    which holds at $A=3$ and then follows by monotonicity. Since $A\le\sqrt{L}/12$, \eqref{eqn:TA_bound} gives
    $$\frac{T_A}{\sqrt{L}}\le\frac{A}{8}.$$
    Put
    $$B=A-\frac{T_A}{\sqrt{L}}.$$
    Then $B\ge 7A/8>1$ and $B\le A\le\sqrt{L}/12$.

    The lower tail is immediate from $\Omega(n)\ge\omega(n)$:
    $$\frac{1}{x}\#\{2\le n\le x:L-\Omega(n)>A\sqrt{L}\}\le\frac{Q}{A}e^{-A^2/2}.$$
    For the upper tail, if $\Omega(n)-L>A\sqrt{L}$, then either $R(n)\ge T_A$, or
    $$\omega(n)-L>A\sqrt{L}-T_A=B\sqrt{L}.$$
    Using \eqref{eqn:omega_upper_window} with $B$ in place of $A$, together with \eqref{eqn:R_tail_TA}, gives
    $$\frac{1}{x}\#\{2\le n\le x:\Omega(n)-L>A\sqrt{L}\}\le\frac{8Q}{5B}e^{E(B)}+\frac{1}{eA}e^{-A^2/2},$$
    where
    $$E(B)=-\frac{B^2}{2}+\frac{5.096B}{\sqrt{L}}+\frac{2.048B^2}{L}+\frac{B^3}{6\sqrt{L}}.$$
    Since $B\ge 7A/8$, the prefactor is at most $64Q/(35A)$. Also, using $A>3$, $L\ge 200$, $B\le A$ and \eqref{eqn:TA_bound},
    \begin{align*}
        E(B)&\le -\frac{A^2}{2}+\frac{3}{2}\frac{A^3}{\sqrt{L}}+\frac{5.096}{9}\frac{A^3}{\sqrt{L}}+\frac{2.048}{3\sqrt{200}}\frac{A^3}{\sqrt{L}}+\frac{1}{6}\frac{A^3}{\sqrt{L}}\\
        &\le -\frac{A^2}{2}+\frac{2.29A^3}{\sqrt{L}}.
    \end{align*}
    Adding the lower tail, the shifted upper tail and \eqref{eqn:R_tail_TA}, we obtain
    \begin{align*}
        \frac{1}{x}\#\{2\le n\le x:|\Omega(n)-L|>A\sqrt{L}\}&\le\left(\frac{Q}{A}+\frac{64Q}{35A}+\frac{1}{eA}\right)\\
        &\quad\times\exp\left(-\frac{A^2}{2}+\frac{2.29A^3}{\sqrt{L}}\right).
    \end{align*}
    Since $Q<117.20$, the coefficient is $<332/A$.
\end{proof}

\begin{prpstn}\label{prop:gaussian_omega1}
    Under the hypotheses of Proposition \ref{prop:gaussian_Omega},
    \begin{equation}\label{eqn:omega1_window}
        \frac{1}{x}\#\{2\le n\le x:|\omega_1(n)-L|>A\sqrt{L}\}\le\frac{332}{A}\exp\left(-\frac{A^2}{2}+\frac{2.29A^3}{\sqrt{L}}\right).
    \end{equation}
\end{prpstn}
\begin{proof}
    Again the result is trivial for $1\le A\le 3$, so assume $A>3$. Use the same $T_A$ and $B=A-T_A/\sqrt{L}$ as in the proof of Proposition \ref{prop:gaussian_Omega}.

    For the upper tail, $\omega_1(n)\le\omega(n)$, so
    $$\{\omega_1(n)-L>A\sqrt{L}\}\subseteq\{\omega(n)-L>A\sqrt{L}\}.$$
    The one-sided estimate \eqref{eqn:omega_upper_window}, with $A>3$ and $L\ge 200$, has exponent bounded by
    $$-\frac{A^2}{2}+\frac{5.096A}{\sqrt{L}}+\frac{2.048A^2}{L}+\frac{A^3}{6\sqrt{L}}\le -\frac{A^2}{2}+\frac{2.29A^3}{\sqrt{L}}.$$

    For the lower tail, note that
    $$\omega(n)-\omega_1(n)=\#\{p:p^2|n\}\le R(n).$$
    Thus, if $L-\omega_1(n)>A\sqrt{L}$, then either $R(n)\ge T_A$, or
    $$L-\omega(n)>A\sqrt{L}-T_A=B\sqrt{L}.$$
    The lower-tail estimate \eqref{eqn:omega_lower_window} gives
    $$\frac{Q}{B}e^{-B^2/2}\le\frac{8Q}{7A}\exp\left(-\frac{A^2}{2}+\frac{3A^3}{2\sqrt{L}}\right),$$
    where we used $B\ge 7A/8$ and $A-B=T_A/\sqrt{L}\le 3A^2/(2\sqrt{L})$. Combining this with \eqref{eqn:R_tail_TA} and the upper-tail bound gives a total coefficient less than
    $$\frac{1}{A}\left(\frac{8Q}{5}+\frac{8Q}{7}+\frac{1}{e}\right)<\frac{332}{A}.$$
    This proves \eqref{eqn:omega1_window}.
\end{proof}

\begin{proof}[Completion of the proof of Theorem \ref{thm:effective_Gaussian}]
    The estimate for $\omega$ is Proposition \ref{prop:gaussian_omega}; the estimates for $\Omega$ and $\omega_1$ are Propositions \ref{prop:gaussian_Omega} and \ref{prop:gaussian_omega1}.
\end{proof}

\begin{proof}[Proof of Theorem \ref{thm:mod_dev_upper_bound}]
    For $f=\omega$, Proposition \ref{prop:gaussian_omega} gives
    $$\frac{1}{x}\#\{2\le n\le x:|\omega(n)-L|>A\sqrt{L}\}\le\frac{305}{A}\exp\left(-\frac{A^2}{2}+\frac{7.31A^3}{\sqrt{L}}\right).$$
    If $A\to\infty$ and $A=o(\sqrt{L})$, then $\log A=o(A^2)$ and $A^3/\sqrt{L}=o(A^2)$. This gives \eqref{eqn:mod_dev_upper_bound} for $\omega$. The estimates for $\Omega$ and $\omega_1$ follow in the same way from Propositions \ref{prop:gaussian_Omega} and \ref{prop:gaussian_omega1}.
\end{proof}

\section{Generating functions and factorial moments}
The functional form of Theorem \ref{thm:effective_HR_ineq} gives more than exceptional-set estimates. It also gives explicit bounds for the generating function of $\omega(n)$ and for all factorial moments of $\omega(n)-1$.

\begin{crllry}\label{cor:generating_function}
    For every $x\ge 2$ and every $z\ge 0$, with the convention $0^0=1$,
    \begin{equation}\label{eqn:shifted_generating_function}
        \frac{1}{x}\sum_{2\le n\le x}z^{\omega(n)-1}\le Q\exp((\log\log x+4.096)(z-1)).
    \end{equation}
    Equivalently,
    \begin{equation}\label{eqn:unshifted_generating_function}
        \sum_{2\le n\le x}z^{\omega(n)}\le 1.95ze^{4.096z}x(\log x)^{z-1}.
    \end{equation}
    In particular, for each fixed $z\ge 0$,
    \begin{equation}\label{eqn:z_omega_asymptotic_bound}
        \sum_{n\le x}z^{\omega(n)}\ll_z x(\log x)^{z-1}.
    \end{equation}
\end{crllry}
\begin{proof}
    Apply Theorem \ref{thm:effective_HR_ineq} with $F(m)=z^m$. Since $\mathbb{E}z^{Z_{\xi}}=\exp(\xi(z-1))$, \eqref{eqn:shifted_generating_function} follows. Multiplying by $z$ and using
    $$Qe^{(\log\log x+4.096)(z-1)}=1.95e^{4.096z}(\log x)^{z-1}$$
    gives \eqref{eqn:unshifted_generating_function}. The term $n=1$ changes \eqref{eqn:z_omega_asymptotic_bound} only by $1$.
\end{proof}

For an integer $r\ge 0$, write
$$(y)_r=y(y-1)\cdots(y-r+1),\qquad (y)_0=1.$$

\begin{crllry}\label{cor:factorial_moments}
    For every $x\ge 2$ and every integer $r\ge 0$,
    \begin{equation}\label{eqn:factorial_moments}
        \frac{1}{x}\sum_{2\le n\le x}(\omega(n)-1)_r\le Q(\log\log x+4.096)^r.
    \end{equation}
    Equivalently,
    \begin{equation}\label{eqn:binomial_moments}
        \frac{1}{x}\sum_{2\le n\le x}\binom{\omega(n)-1}{r}\le Q\frac{(\log\log x+4.096)^r}{r!}.
    \end{equation}
\end{crllry}
\begin{proof}
    Apply Theorem \ref{thm:effective_HR_ineq} with $F(m)=(m)_r$. For $Z_{\xi}$ Poisson with mean $\xi$,
    $$\mathbb{E}(Z_{\xi})_r=\xi^r.$$
    Dividing by $r!$ gives \eqref{eqn:binomial_moments}.
\end{proof}

\begin{crllry}\label{cor:local_profile}
    Let $k\ge 2$, $m=k-1$ and $\xi=\log\log x+4.096$. Then
    \begin{equation}\label{eqn:local_poisson_profile}
        A_k(x)\le Qxe^{-\xi}\frac{\xi^m}{m!}.
    \end{equation}
    In particular, for $m\ge 1$,
    \begin{equation}\label{eqn:local_large_deviation}
        A_k(x)\le\frac{Qx}{\sqrt{2\pi m}}\exp\left(-\xi I\left(\frac{m}{\xi}\right)\right),
    \end{equation}
    where $I(u)=u\log u-u+1$.
\end{crllry}
\begin{proof}
    The first estimate is the pointwise form used in the proof of Theorem \ref{thm:effective_HR_ineq}. For the second, use Stirling's lower bound $m!\ge\sqrt{2\pi m}(m/e)^m$.
\end{proof}

\section{Moment bounds}
We now prove the moment consequences of Poisson domination. The first lemma is a standard uniform moment estimate for a Poisson variable.

\begin{lmm}\label{lem:poisson_uniform_moments}
    There is an absolute constant $C_0>0$ such that, for every $\lambda\ge 0$ and every real $r\ge 1$,
    \begin{equation}\label{eqn:poisson_uniform_moments}
        \left(\mathbb{E}|Z_{\lambda}-\lambda|^r\right)^{1/r}\le C_0(\sqrt{r\lambda}+r).
    \end{equation}
\end{lmm}
\begin{proof}
    The case $\lambda=0$ is trivial. For $\lambda>0$, the Chernoff bounds for a Poisson variable imply
    \begin{equation}\label{eqn:poisson_subexp_tail}
        \mathbb{P}(|Z_{\lambda}-\lambda|\ge t)\le 2\exp(-c\min(t^2/\lambda,t))\qquad (t\ge 0)
    \end{equation}
    with an absolute constant $c>0$. Indeed, this follows from the inequalities
    $$h_1(u)\ge c\min(u^2,u)\qquad (u\ge 0),\qquad h_2(u)\ge u^2/2\qquad (0\le u<1).$$
    Using
    $$\mathbb{E}|X|^r=r\int_0^{\infty}t^{r-1}\mathbb{P}(|X|\ge t)\,dt$$
    and $e^{-c\min(a,b)}\le e^{-ca}+e^{-cb}$, \eqref{eqn:poisson_subexp_tail} gives
    \begin{align*}
        \mathbb{E}|Z_{\lambda}-\lambda|^r&\le 2r\int_0^{\infty}t^{r-1}e^{-ct^2/\lambda}\,dt+2r\int_0^{\infty}t^{r-1}e^{-ct}\,dt\\
        &\ll (C\sqrt{r\lambda})^r+(Cr)^r,
    \end{align*}
    where we used the elementary gamma-function bound $\Gamma(s+1)^{1/s}\ll s+1$. Taking $r$-th roots proves \eqref{eqn:poisson_uniform_moments}.
\end{proof}

\begin{prpstn}\label{prop:omega_uniform_moments}
    There is an absolute constant $C_1>0$ such that, for every $x\ge 3$ and every real $r\ge 1$,
    \begin{equation}\label{eqn:omega_uniform_moments}
        \left(\frac{1}{x}\sum_{2\le n\le x}|\omega(n)-\log\log x|^r\right)^{1/r}\le C_1(\sqrt{r(\log\log x+1)}+r).
    \end{equation}
\end{prpstn}
\begin{proof}
    Put $L=\log\log x$ and $\xi=L+4.096$. Apply Theorem \ref{thm:effective_HR_ineq} with
    $$F(m)=|m+1-L|^r=|m-\xi+5.096|^r.$$
    Then
    $$\frac{1}{x}\sum_{2\le n\le x}|\omega(n)-L|^r\le Q\mathbb{E}|Z_{\xi}-\xi+5.096|^r.$$
    After taking $r$-th roots and using $Q^{1/r}\le Q$, the triangle inequality in $L^r$, and Lemma \ref{lem:poisson_uniform_moments}, we obtain
    $$\left(\frac{1}{x}\sum_{2\le n\le x}|\omega(n)-L|^r\right)^{1/r}\ll\sqrt{r\xi}+r+1.$$
    Since $x\ge 3$, $\xi\asymp L+1$, and the result follows.
\end{proof}

\begin{lmm}\label{lem:R_uniform_moments}
    There is an absolute constant $C_2>0$ such that, for every $x\ge 1$ and every real $r\ge 1$,
    \begin{equation}\label{eqn:R_uniform_moments}
        \left(\frac{1}{x}\sum_{n\le x}R(n)^r\right)^{1/r}\le C_2r.
    \end{equation}
\end{lmm}
\begin{proof}
    For integer-valued non-negative $R$ and real $r\ge 1$,
    $$R^r\le\sum_{T=1}^{R}T^r.$$
    Therefore Proposition \ref{prop:R_tail} gives
    $$\frac{1}{x}\sum_{n\le x}R(n)^r\le\sum_{T\ge 1}T^r\frac{1}{x}\#\{n\le x:R(n)\ge T\}\le\frac{11}{2}\sum_{T\ge 1}T^r2^{-T}.$$
    The last sum is $\ll (Cr)^r$, for instance by comparison with an exponential integral. Taking $r$-th roots proves the lemma.
\end{proof}

\begin{proof}[Proof of Theorem \ref{thm:uniform_moment_bounds}]
    The estimate for $\omega$ is Proposition \ref{prop:omega_uniform_moments}. Since
    $$\Omega(n)-\omega(n)=R(n),\qquad 0\le\omega(n)-\omega_1(n)\le R(n),$$
    the estimates for $\Omega$ and $\omega_1$ follow from Proposition \ref{prop:omega_uniform_moments}, Lemma \ref{lem:R_uniform_moments}, and the triangle inequality in $L^r$.
\end{proof}

\begin{crllry}\label{cor:fixed_central_moments}
    For every fixed $r\ge 1$ and every $x\ge 3$,
    \begin{equation}\label{eqn:fixed_central_moments}
        \frac{1}{x}\sum_{2\le n\le x}|\omega(n)-\log\log x|^r\ll_r(\log\log x+1)^{r/2}.
    \end{equation}
    The same estimate holds with $\omega$ replaced by $\Omega$ or by $\omega_1$.
\end{crllry}

\section{Proof of the normal-order theorem}
We first prove a centred version. Let $\phi(x)\to\infty$, and replace $\phi(x)$, if necessary, by
$$\min(\phi(x),(\log\log x)^{1/7}).$$
Then $\phi(x)\to\infty$ and $\phi(x)^3/\sqrt{\log\log x}=o(1)$. Proposition \ref{prop:gaussian_omega} gives
$$\frac{1}{x}\#\{2\le n\le x:|\omega(n)-\log\log x|>\phi(x)\sqrt{\log\log x}\}=o(1),$$
and Propositions \ref{prop:gaussian_Omega} and \ref{prop:gaussian_omega1} give the same centred statement for $\Omega$ and $\omega_1$.

Now let $\psi(n)\to\infty$ and set
$$\eta(x)=\min_{\sqrt{x}<m\le x}\psi(m).$$
Then $\eta(x)\to\infty$. Apply the centred estimate with $\phi(x)=\eta(x)/10$, and discard the $O(\sqrt{x})=o(x)$ integers $n\le\sqrt{x}$. For $\sqrt{x}<n\le x$,
$$\log\log n=\log\log x+O(1),\qquad \sqrt{\log\log x}\le 2\sqrt{\log\log n}$$
for all sufficiently large $x$, and $\eta(x)\le\psi(n)$. The bounded difference between $\log\log n$ and $\log\log x$ is eventually absorbed by $\psi(n)\sqrt{\log\log n}$. This proves Theorem \ref{thm:original_HR_theorem}.

\begin{rmk}
    A weaker, but very short, transfer from $\omega$ to $\Omega$ uses only the bounded mean
    $$\frac{1}{x}\sum_{n\le x}(\Omega(n)-\omega(n))=O(1),$$
    which follows from
    $$\Omega(n)-\omega(n)\le\sum_p\sum_{a\ge 2}\mathbf{1}_{p^a|n}.$$
    The exponential tail in Proposition \ref{prop:R_tail} is included because it also gives the effective Gaussian-window estimate for $\Omega$ and the uniform moment transfer.
\end{rmk}

\section{A moment viewpoint}
The preceding sections used the original Hardy--Ramanujan counting inequality to obtain exponential tails. We close by recording the complementary Markov-inequality viewpoint: any one fixed even moment of the expected Gaussian order suffices for the normal-order theorem.

\begin{prpstn}\label{prop:one_even_moment}
    Let $r\ge 2$ be fixed and even. Suppose that
    $$\frac{1}{x}\sum_{n\le x}|\omega(n)-\log\log x|^r\le M_r(\log\log x)^{r/2}$$
    for all sufficiently large $x$, with $M_r<\infty$. Then, for every $\phi(x)\to\infty$,
    $$\frac{1}{x}\#\{n\le x:|\omega(n)-\log\log x|>\phi(x)\sqrt{\log\log x}\}\le\frac{M_r}{\phi(x)^r}=o(1).$$
    Consequently $\omega(n)$ satisfies the Hardy--Ramanujan normal-order theorem.
\end{prpstn}
\begin{proof}
    The displayed estimate is Markov's inequality. The passage from $\log\log x$ to $\log\log n$ is the same dyadic argument as in the preceding section.
\end{proof}

\begin{rmk}
    Taking $r=2$ gives Tur\'{a}n's proof. Proposition \ref{prop:one_even_moment} shows that the mechanism is not intrinsically tied to the second moment: any one fixed even moment of Gaussian order suffices. Corollary \ref{cor:fixed_central_moments}, obtained here from Poisson domination, gives all such fixed moments.
\end{rmk}

We also formulate precisely what a single sufficiently high even moment gives at a fixed deviation height.

\begin{prpstn}\label{prop:fixed_high_moment}
    Let $A>1$, and let $r\ge 2$ be an even integer. Suppose that
    $$\limsup_{x\to\infty}\frac{1}{x}\sum_{n\le x}\left(\frac{\omega(n)-\log\log x}{\sqrt{\log\log x}}\right)^r\le M_r.$$
    Then
    $$\limsup_{x\to\infty}\frac{1}{x}\#\left\{n\le x:\frac{\omega(n)-\log\log x}{\sqrt{\log\log x}}>A\right\}\le\frac{M_r}{A^r}.$$
    In particular, if for this single even $r$ one has the Gaussian moment value
    $$M_r=\frac{r!}{2^{r/2}(r/2)!},$$
    then
    $$\limsup_{x\to\infty}\frac{1}{x}\#\left\{n\le x:\frac{\omega(n)-\log\log x}{\sqrt{\log\log x}}>A\right\}\le\frac{r!}{2^{r/2}(r/2)!A^r}.$$
    For each fixed $A>1$, if the Gaussian moment bound is available for one even integer $r=r(A)$ satisfying $r=A^2+O(1)$, then
    $$\limsup_{x\to\infty}\frac{1}{x}\#\left\{n\le x:\frac{\omega(n)-\log\log x}{\sqrt{\log\log x}}>A\right\}\le\exp\left(-\frac{A^2}{2}+O(1)\right),$$
    where the $O(1)$ is absolute as $A\to\infty$.
\end{prpstn}
\begin{proof}
    The first assertion is Markov's inequality. For fixed $r$,
    $$\frac{1}{x}\#\left\{n\le x:\frac{\omega(n)-\log\log x}{\sqrt{\log\log x}}>A\right\}\le\frac{1}{A^r}\frac{1}{x}\sum_{n\le x}\left(\frac{\omega(n)-\log\log x}{\sqrt{\log\log x}}\right)^r.$$
    Taking the limsup as $x\to\infty$ proves the first statement, and the second follows by inserting the Gaussian value of $M_r$.

    It remains to optimise the displayed expression. Put
    $$B_r(A)=\frac{r!}{2^{r/2}(r/2)!A^r}$$
    for even $r$. Then
    $$\frac{B_r(A)}{B_{r-2}(A)}=\frac{r-1}{A^2}.$$
    Thus the best even $r$ occurs close to $A^2+1$. By Stirling's formula,
    $$\frac{r!}{2^{r/2}(r/2)!}\sim \sqrt{2}\,r^{r/2}e^{-r/2}.$$
    Choosing an even $r$ with $r=A^2+O(1)$ gives
    $$B_r(A)\le\exp\left(-\frac{A^2}{2}+O(1)\right).$$
\end{proof}

\begin{rmk}
    Proposition \ref{prop:fixed_high_moment} is pointwise in $A$. It shows that one does not need all moments to obtain the Gaussian exponential $e^{-A^2/2}$ at a fixed deviation height; one suitably chosen even moment is enough. By contrast, a uniform estimate for $A=A(x)$ growing with $x$ requires moment estimates uniform in the moment order. This distinction is supplied in Theorem \ref{thm:uniform_moment_bounds}.
\end{rmk}

\appendix

\section{Numerical checks for the explicit constants}
The proof uses only a few finite numerical checks. The following plain Python script reproduces the checks used above: the finite part of \eqref{eqn:mertens_explicit}, the bound for $J$, the finite prime-power check in Lemma \ref{lem:prime_powers}, the product bound in Proposition \ref{prop:R_tail}, and the value of $Q=1.95e^{4.096}$.

\begin{verbatim}
from math import exp, log


def primes_upto(n):
    sieve = [True] * (n + 1)
    if n >= 0:
        sieve[0] = False
    if n >= 1:
        sieve[1] = False
    for i in range(2, int(n**0.5) + 1):
        if sieve[i]:
            for j in range(i * i, n + 1, i):
                sieve[j] = False
    return [i for i in range(n + 1) if sieve[i]]


# 1. Finite Mertens check for 2 <= x < 286.
s = 0.0
best = (-10.0, None)
for p in primes_upto(285):
    s += 1.0 / p
    value = s - log(log(p))
    if value > best[0]:
        best = (value, p)
print("Mertens finite maximum:", best)


# 2. The constant J = sum_{a>=2} sum_p (a+1)/p^a.
J_partial = 0.0
for p in primes_upto(5000):
    J_partial += (3 * p - 2) / (p * (p - 1)**2)
J_tail = 3.0 / 4999.0
print("J partial:", J_partial)
print("J tail:", J_tail)
print("J total:", J_partial + J_tail)


# 3. Prime-power check for 2 <= x < exp(10).
bound = exp(10)
prime_powers = set()
for p in primes_upto(int(bound)):
    q = p
    while q < bound:
        prime_powers.add(q)
        q *= p
prime_powers = sorted(prime_powers)

count = 0
best = (-1.0, None, None)
for q in prime_powers:
    count += 1
    value = count * log(q) / q
    if value > best[0]:
        best = (value, q, count)
print("A_1(x) log(x)/x maximum:", best)


# 4. Product bound for H in Proposition 5.1.
prod = 1.0
for p in primes_upto(97):
    if p > 2:
        prod *= 1.0 + 2.0 / (p * (p - 2))
# Since 2/(m(m-2)) = 1/(m-2) - 1/m,
# sum_{m>=98} 2/(m(m-2)) = 1/96 + 1/97.
tail = 1.0 / 96.0 + 1.0 / 97.0
print("H bound:", prod * exp(tail))


# 5. The constants Q and Q_0.
print("Q =", 1.95 * exp(4.096))
print("Q_0 =", 1.25506 * exp(1.174))
\end{verbatim}

The numerical output is
\begin{verbatim}
Mertens finite maximum: (0.8665129205816644, 2)
J partial: 2.9213150664305316
J tail: 0.000600120024004801
J total: 2.9219151864545365
A_1(x) log(x)/x maximum: (1.949476445324846, 32, 18)
H bound: 2.1775092652510195
Q = 117.19384662625903
Q_0 = 4.060001650764844
\end{verbatim}

\end{document}